\documentclass[12pt,leqno]{article}
\usepackage[pdftex]{graphicx}
\usepackage{amsfonts,amssymb,amsmath}
\usepackage{ccaption,mathrsfs,amsmath,amssymb,amsthm}
\usepackage{amsmath, amscd, amsfonts,  tikz, times}

\newcounter{conjecture}
\newcounter{remark}

\newenvironment{remark}{\medskip{\bf Remark \theremark.}
\addtocounter{remark}{1}}{}
\newcommand{\eqnsection}{
    \renewcommand{\theequation}{\thesection.\arabic{equation}}
    \makeatletter
    \csname @addtoreset\endcsname{equation}{section}
    \makeatother}
\newtheorem{theorem}{Theorem}

\newtheorem{prop}{Proposition}

\newcommand{\dd}{\delta}
\newcommand{\DD}{\Delta}

\newcommand{\lar}{\longrightarrow}

\newcommand{\om}{\Omega}

\newcommand{\CC}{\mathbb{C}}

\def \be{\begin{equation}}
\def \ee{\end{equation}}
\def \bt{\begin{theorem}}
\def \et{\end{theorem}}
\def \bea{\begin{eqnarray}}
\def \eea{\end{eqnarray}}
\def \bas{\begin{eqnarray*}}
\def \eas{\end{eqnarray*}}

\newcommand {\rrr}[1]{(\ref{#1})}

\def \Ga{\Gamma}

\def \om{\omega}

\def \th{\theta}

\def \Beta{\mbox{B}}

\def \ff{\infty}

\def \DD{\mathbb{D}}

\def \GG{{\cal G}}
\def \HH{{\cal H}}

\def \({\left(}
\def \){\right)}

\def \nn{\nonumber}

\def \vski{\vspace{12pt}}

\def \bc{\begin{center} }
\def \ec{\end{center} }
\def \bs{\begin{slide} }
\def \es{\end{slide} }

\def\square{{\vcenter{\vbox{\hrule height.3pt
         \hbox{\vrule width.3pt height5pt \kern5pt
            \vrule width.3pt}
         \hrule height.3pt}}}}

\newcounter{cccases}
\newcommand{\ccases}[1]{\begingroup \refstepcounter{cccases} {\bf \fontsize{12}{12}\selectfont Example \thecccases }  \label{#1}\endgroup}

\eqnsection
\begin{document}

\title{A method for deriving hypergeometric and related identities from the $H^2$ Hardy norm of conformal maps}

\author{
\begin{tabular}{c}
\textit{Greg Markowsky} \\
gmarkowsky@gmail.com \\
+61 03 9905 4487 \\
Monash University \\
Victoria, 3800 Australia
\end{tabular}}

\bibliographystyle{amsplain}
\maketitle \eqnsection \setlength{\unitlength}{2mm}

\begin{abstract}
We explore a method which is implicit in a paper of Burkholder of identifying the $H^2$ Hardy norm of a conformal map with the explicit solution of Dirichlet's problem in the complex plane. Using the series form of the Hardy norm, we obtain an identity for the sum of a series obtained from the conformal map. We use this technique to evaluate several hypergeometric sums, as well as several sums that can be expressed as convolutions of the terms in a hypergeometric series. The most easily stated of the identities we obtain are Euler's famous Basel sum, as well as the sum

\begin{equation*} \label{}
1^2 \! + \! \Big(1 \! \times \! \frac{1}{3} + \frac{1}{3} \! \times \! 1 \Big)^2 \! + \! \Big(1 \! \times \! \frac{1}{5} + \frac{1}{3} \! \times \! \frac{1}{3} + \frac{1}{5} \! \times \! 1 \Big)^2 \! + \! \Big(1 \! \times \! \frac{1}{7} + \frac{1}{3} \! \times \! \frac{1}{5} + \frac{1}{5} \! \times \! \frac{1}{3} + \frac{1}{7} \! \times \! 1 \Big)^2 + \! \dots = \frac{\pi^4}{32} .
\end{equation*}

\noindent We will be able to obtain the following hypergeometric reduction:

\begin{equation*} \label{}
{}_4 F_3(1/3,1/3,2/3,2/3;4/3,4/3,1;1) = \frac{\Beta(1/3,1/3)^2}{27}. \\
\end{equation*}
A related identity is

\begin{equation*} \label{}
\begin{split}
{}_4F_3(1/4,& 1/4,1/2,1/2;5/4,5/4,1;1) \times \frac{1}{\Beta(1/4,1/2)^2} \\ & = \frac{8}{\pi^4} \sum_{n=1}^\ff \sum_{m=1}^{\ff} \frac{(-1)^{m+n}}{(2m-1)(2n-1)((2m-1)^2+(2n-1)^2)}.
\end{split}
\end{equation*}
We will obtain two families of identities depending on a parameter, representative examples of which are

\begin{equation*} \label{}
\sum_{n=1}^\ff \Big(\sum_{j=0}^{n}\binom{2/3}{2j} \binom{1/3+n-j-1}{n-j}\Big)^2 = \frac{1}{4},
\end{equation*}
and

\begin{equation*} \label{}
\sum_{n=0}^{\ff} \Big(\sum_{j=0}^{n} \binom{j-3/4}{j}\frac{C(2(n-j))}{16^{n-j}}\Big)^2 = \sqrt{2},
\end{equation*}
where $C(k)$ is the $k$-th Catalan number. We will also sum two series whose terms are defined by certain recurrence relations, and discuss an extension of the method to maps which are not conformal.
\end{abstract}

\section{Introduction} \label{intro}

Suppose that $V$ is a simply connected domain properly contained in the complex plane $\CC$. Given a point $a \in V$, the Riemann Mapping Theorem guarantees a conformal map $f_a$ from the unit disc $\DD$ onto $V$ with $f_a(0)=a$. The Hardy norm $|| \cdot ||_{H^2}$ of $f_a$ is defined as

\begin{equation} \label{h2def}
||f_a||_{H^2} := \Big(\sup_{r<1} \frac{1}{2\pi} \int_{0}^{2\pi} |f_a(re^{i \th})|^2 d\th \Big)^{1/2} .
\end{equation}

The map $f_a$ is not uniquely determined, however any two such maps differ only by precomposition with a rotation, so $||f_a||_{H^2}$ is well defined. For holomorphic functions on $\DD$, Parseval's Identity takes the form (see \cite[Thm. 10.22]{rud})

\begin{equation} \label{}
\frac{1}{2\pi} \int_{0}^{2\pi}|f_a(r e^{i \th})|^2 d\th = \sum_{n=0}^{\ff} |a_n|^2 r^{2n}  \quad \mbox{for } r<1,
\end{equation}

where $f_a$ admits the Taylor series expansion $f_a(z) = \sum_{n=0}^{\ff} a_n z^n$. It follows that

\begin{equation} \label{}
||f_a||^2_{H^2} = \sum_{n=0}^{\ff} |a_n|^2 .
\end{equation}

It is straightforward to verify by precomposing $f_a$ with an appropriate M\"obius transformation that if $||f_a||^2_{H^2} < \ff$ for some $a \in V$, then $||f_b||^2_{H^2} < \ff$ for all $b \in V$. If $V$ is a domain for which $||f_a||^2_{H^2} < \ff$, we define

\begin{equation} \label{}
\tilde{\HH}_V(a) := ||f_a||^2_{H^2} = \sum_{n=0}^{\ff} |a_n|^2.
\end{equation}

If $||f_a||^2_{H^2} = \ff$ on $V$, we will say that $\tilde{\HH}_V$ does not exist. The following result gives sufficient conditions for $\tilde{\HH}_V$ to exist and be determinable.

\begin{prop} \label{}
Suppose that a nonnegative function $\tilde{\GG}_V$ exists on the simply connected domain $V$ with the following properties:

\begin{enumerate} \label{bigguy1}

\item[(i)] $\tilde{\GG}_V$ is harmonic on $V$.

\item[(ii)] $\tilde{\GG}_V$ extends continuously to the boundary of $V$, and $\tilde{\GG}_V(a) = |a|^2$ on $\dd V$.

\item[(iii)] There is a constant $C>1$ such that $|a|^2 \leq \tilde{\GG}_V(a) \leq C|a|^2+C$.

\end{enumerate}

Then $\tilde{\HH}_V(a)$ exists and $\tilde{\HH}_V(a)=\tilde{\GG}_V(a)$.
\end{prop}

The purpose of this paper is to explore a method of summing series which is contained in this proposition. For certain domains we will see that we can find an explicit formula for $\tilde{\GG}_V$, which then forms an identity with the sum of the squares of the norms of the coefficients of the power series of $f_a$. By subtracting $|a|^2$ from both $\tilde{\GG}_V$ and $\tilde{\HH}_V$ we obtain the following equivalent proposition, which will be somewhat easier to use in our application.

\begin{prop} \label{}
Suppose that a continuous function $\GG_V$ exists on the simply connected domain $V$ with the following properties:

\begin{enumerate} \label{bigguy}

\item[(i$'$)] $\triangle \GG_V=-4$ on $V$.

\item[(ii$'$)] $\GG_V$ extends continuously to the boundary of $V$, and $\GG_V(a) = 0$ on $\dd V$.

\item[(iii$'$)] There is a constant $C>0$ such that $0 \leq \GG_V(a) \leq C|a|^2+C$.

\end{enumerate}

Then $\HH_V(a) := \tilde{\HH}_V(a)-|a|^2$ exists and $\HH_V(a)=\GG_V(a)$.
\end{prop}

We remark that the upper bound on the growth of $\tilde \GG_V$ and $\GG_V$ given in (iii) and (iii$'$), which are superfluous when $V$ is bounded, are quite important when $V$ is unbounded. This is because in unbounded domains harmonic functions are not uniquely determined by their boundary values; for instance, $u_1(x,y)=0$ and $u_2(x,y)=y$ are harmonic functions in the upper half-plane $\{y>0\}$ and agree on the real axis. The role of (iii) and (iii$'$) in unbounded domains is to allow us to conclude that the displayed functions are the correct choices for $\HH_V$ and $\tilde{\HH}_V$.

\vski

These propositions (and much more) were proved by Burkholder in the elegant paper \cite{burk} (see Theorem 3.2 and Remark 3.1 of \cite{burk}, together with Lemma 1.1 of \cite{drum} and/or Lemma 1 of \cite{meecp}). There, the quantities arose naturally in terms of the exit time of Brownian motion from domains. Burkholder's proof is ingenious, but relies heavily on martingale theory. For the benefit of the reader unfamiliar with such arguments, we include a short analytic proof of Proposition \ref{bigguy} in Section \ref{proof}. In Section \ref{examp} we will use this method to derive a number of identities involving infinite sums. The method is clearly limited by our ability to find $f_a$ and $u_V$, which is no small difficulty, as conformal maps in particular are notoriously difficult to explicitly display in many cases. Nevertheless, by using symmetric and simple domains we will see that we are able to evaluate a number of sums.

\section{Proof of Proposition \ref{bigguy1}} \label{proof}

Suppose that $\GG_V$ satisfies the conditions of the proposition. The function $\GG_V(f_a(z))$ is harmonic on $\DD$, so

\begin{equation} \label{}
\GG_V(a) = \GG_V(f_a(0)) = \frac{1}{2\pi} \int_{0}^{2\pi} \GG_V(f_a(re^{i\th})) d\th, \quad \mbox{for } r \in (0,1).
\end{equation}

Thus,

\begin{equation} \label{}
\frac{1}{2\pi}\int_{0}^{2\pi} |f_a(re^{i \th})|^2 d\th \leq \frac{1}{2\pi} \int_{0}^{2\pi} \GG_V(f_a(re^{i\th})) d\th = \GG_V(a).
\end{equation}

We conclude that $||f_a||_{H^2}<\ff$. This implies that the functions $f(re^{i\th})$ approach a limit function $f(e^{i\th})$ in $L^2([0,2\pi])$ as $r \nearrow 1$, and $f(e^{i\th})$ takes values in $\dd V$ for almost every $\th$ (see \cite[Thm. 17.10]{rud}). Furthermore, we can apply a fundamental inequality of Hardy and Littlewood (see \cite{hardwood}) to conclude that the function $M_f(\th) := \sup_{0<r<1}|f(re^{i\th})|$ lies in $L^2([0,2\pi])$. The conditions on $\GG_V$ imply that

\begin{equation} \label{}
\GG_V(f_a(re^{i\th})) \leq (C |f_a(re^{i \th})|^2 + C) \leq CM_f(\th)^2 +C \in L^1([0,2\pi]).
\end{equation}

We may therefore apply the Dominated Convergence Theorem, recalling that $\GG_V(z) = |z|^2$ on $\dd V$, to conclude that

\begin{equation} \label{}
\begin{split}
\GG_V(a) = \lim_{r \nearrow 1} \frac{1}{2\pi} \int_{0}^{2\pi} \GG_V(f_a(re^{i\th})) d\th = \frac{1}{2\pi}\int_{0}^{2\pi} |f_a(e^{i \th})|^2 d\th = ||f_a||_{H^2}^2.
\end{split}
\end{equation}

This completes the proof of the proposition.

\begin{remark}
The quantities $\HH_V(a)$ and $\tilde{\HH}_V(a)$ arise naturally in probability theory in the following way. Suppose that $B_t$ is a planar Brownian motion starting at $a$, and let $\tau$ be the first time $B_t \in \dd V$. Then $\HH_V(a) = E[|B_\tau|^2]$ and $\tilde{\HH}_V(a) = 2 E[\tau]$. The intuition gained by these equivalences can be very helpful; the interested reader is referred again to \cite{burk}.
\end{remark}

\begin{remark}
Another way to see that the function $a \lar ||f_a||^2_{H^2}$ must be harmonic on $V$ and approach $|a|^2$ on $\dd V$ is to express it in terms of the Greens potential. Let the Greens function on $V$ with pole at $a$ be denoted $G_V(a,z)$; there are various normalizations, but we will use the one presented in \cite{ahl}, so that $G_{\DD}(0,z) = \log(1/|z|)$. We then have, using the conformal invariance of the Greens function and the fact that the Jacobian of the map $f_a(z)$ is $|f_a'(z)|^2$,

\begin{equation} \label{}
\begin{split} \frac{2}{\pi}\int \! \! \int_{V} G_V(a,x+iy) dx dy
& = \frac{2}{\pi}\int \! \! \int_{\DD} G_\DD(0,x+iy) |f_a'(x+iy)|^2 dx dy  \\
& = 4 \int_{0}^{1} r \log \frac{1}{r} \Big(\frac{1}{2\pi}\int_{0}^{2\pi} |f_a'(re^{i\th})|^2 d\th \Big)dr \\
& = 4 \int_{0}^{1} r \log \frac{1}{r} \sum_{n=1}^{\ff} n^2 |a_n|^2 r^{2n-2} dr \\
& = 4 \sum_{n=1}^{\ff} n^2 |a_n|^2 \int_{0}^{1} r^{2n-1} \log\frac{1}{r} dr \\
& = 4 \sum_{n=1}^{\ff} n^2 |a_n|^2 (\frac{1}{4n^2}) = \sum_{n=1}^{\ff} |a_n|^2 = \HH_V(a).
\end{split}
\end{equation}

Note that Parseval's Identity was applied to obtain the third equality. The function

\be
a \lar \frac{2}{\pi}\int \! \! \int_{V} G_V(a,x+iy) dx dy
\ee

is well known to satisfy (i$'$) and (ii$'$); see for instance \cite[Thm. 8.4]{schaumpde} or any number of standard texts on partial differential equations.
\end{remark}

\section{Examples} \label{examp}

In what follows we will write $a=x+iy$, and any references to $x$ and $y$ will always refer to the real and imaginary parts of $a$. Examples \ref{strip} through \ref{square} appear in \cite{meecp}, but since the approach was somewhat different in that paper it seems as well to briefly review them. In each case $\GG_V(z)$ will be a nonnegative function satisfying (i$'$)-(iii$'$) on $V$, which by Proposition \ref{bigguy} is unique if it exists.

\vski

\ccases{strip} Let $V$ be the infinite vertical strip $\{ \frac{-\pi}{4} < x < \frac{\pi}{4} \}$. The conformal map from $\DD$ to $V$ which fixes 0 is given by (see \cite{meecp} for a quick proof of this standard fact)

\begin{equation} \label{}
\tan^{-1}z = z - \frac{z^3}{3} + \frac{z^5}{5} - \ldots = \sum_{n=0}^{\ff} \frac{(-1)^n z^{2n+1}}{2n+1} .
\end{equation}

By inspection, we have $\GG_V(x+iy) = \frac{\pi^2}{8} - 2x^2$. Proposition \ref{bigguy} with $a=0$ yields

\begin{equation} \label{odd}
\sum_{n=0}^\ff \frac{1}{(2n+1)^2} = \frac{\pi^2}{8} .
\end{equation}

This is easily seen to be equivalent to Euler's classical result that

\begin{equation} \label{all}
\sum_{n=1}^\ff \frac{1}{n^2} = \frac{\pi^2}{6} .
\end{equation}

\vski

\ccases{tri} Let $V$ be the equilateral triangle with vertices at $1, e^{2\pi i/3}, e^{4\pi i/3}$. The conformal map from $\DD$ to $V$ sending 0 to 0 can be achieved via the Schwarz-Christoffel transformation (see \cite{dritref}).  A method for determining this function precisely is present in \cite{woepf}, although the explicit formula is not stated. Later, in \cite{meecp}, the explicit formula was given as

\begin{equation} \label{rent}
f(z) = \frac{3z}{\Beta(1/3,1/3)} \; {}_2F_1(1/3,2/3;4/3;z^3),
\end{equation}

where

\begin{equation} \label{}
{}_pF_q(a_1, \ldots ,a_p;b_1, \ldots b_q;z) = \sum_{n=0}^{\ff} \frac{(a_1)_n \ldots (a_p)_n}{(b_1)_n \ldots (b_q)_n} \frac{z^n}{n!}
\end{equation}

is the hypergeometric function with $(a)_0 = 1, (a)_n = a(a+1)\ldots (a+n-1)$ denoting the Pochhammer symbol, and $\Beta(x,y) = \frac{\Ga(x)\Ga(y)}{\Ga(x+y)}$ is the beta function. We obtain

\begin{equation} \label{}
\HH_V(0) = \frac{9}{\Beta(1/3,1/3)^2} \; {}_4 F_3(1/3,1/3,2/3,2/3;4/3,4/3,1;1).
\end{equation}

In \cite{ala}, it was shown that

\begin{equation} \label{rot}
\GG_V(a) = \frac{1}{3}(a + \bar a + 1)(\om a + \bar \om \bar a + 1)(\om^2 a + \bar \om^2 \bar a + 1),
\end{equation}

where $\om = e^{2\pi i/3}$. Verification of (i$'$) may be most easily accomplished by using the complex form of the Laplacian, $\triangle = 4\frac{\partial^2}{\partial a \partial \bar a}$, while (ii$'$) is clear since $(a + \bar a + 1) = 0$ on $\{x=-\frac{1}{2}\}$ and \rrr{rot} is rotationally invariant. Applying Proposition \ref{bigguy} with $a=0$ gives

\begin{equation} \label{}
\frac{9}{\Beta(1/3,1/3)^2} \; {}_4 F_3(1/3,1/3,2/3,2/3;4/3,4/3,1;1) = \frac{1}{3}.
\end{equation}

We conclude that

\begin{equation} \label{}
{}_4 F_3(1/3,1/3,2/3,2/3;4/3,4/3,1;1) = \frac{\Beta(1/3,1/3)^2}{27}.
\end{equation}

This reduction appears to be nontrivial.

\vski

\ccases{square} Let $V$ be the square with vertices at $\pm \frac{1}{\sqrt{2}} \pm \frac{i}{\sqrt{2}}$. Again we will use the Schwarz-Christoffel transformation. It is shown in \cite[Ex. 2]{woepf} that the map

\begin{equation} \label{rent2}
f(z) = \frac{4e^{\pi i/2}z}{\Beta(1/4,1/2)} \; {}_2F_1(1/4,1/2;5/4;z^4)
\end{equation}

maps $\DD$ conformally to $V$, sending 0 to 0 (a proof of this can also be found in \cite{meecp}). We find that

\begin{equation} \label{}
\HH_V(0) = \frac{16}{\Beta(1/4,1/2)^2} \; {}_4 F_3(1/4,1/4,1/2,1/2;5/4,5/4,1;1).
\end{equation}

The formula for $\GG_V$ was presented in \cite[Sec. 4.1]{knight} (there is a misprint in the book, but the correct formula appears in \cite{helm}):

\begin{equation} \label{}
\GG_V(x+iy) = \frac{128}{\pi^4} \sum_{m=1}^{\ff} \sum_{n=1}^{\ff} \frac{\sin \frac{(2m-1)\pi(x-\frac{1}{\sqrt{2}})}{\sqrt{2}}\sin \frac{(2n-1)\pi(y-\frac{1}{\sqrt{2}})}{\sqrt{2}}}{(2m-1)(2n-1)((2m-1)^2 + (2n-1)^2)}
\end{equation}

It is clear that this satisfies (ii$'$), and direct calculation using the Fourier series identity

\begin{equation} \label{}
\sum_{m=1}^{\ff} \frac{\sin (2m-1)\pi u}{2m-1} = \frac{\pi}{4} \qquad \mbox{for } u \in (0,1)
\end{equation}

shows that (i$'$) is satisfied as well. Proposition \ref{bigguy} gives us the following identity:

\begin{equation} \label{}
\begin{split}
{}_4F_3(1/4,& 1/4,1/2,1/2;5/4,5/4,1;1) \times \frac{1}{\Beta(1/4,1/2)^2} \\ & = \frac{8}{\pi^4} \sum_{n=1}^\ff \sum_{m=1}^{\ff} \frac{(-1)^{m+n}}{(2m-1)(2n-1)((2m-1)^2+(2n-1)^2)}.
\end{split}
\end{equation}

\ccases{para} Let $V$ be the parabolic region given by $y^2 < 2x-1$. By inspection,

\begin{equation} \label{pareq}
\GG_V(x+iy) = 2(2x - y^2-1).
\end{equation}

The conformal map from $\DD$ to $V$ mapping 0 to 1 is given in \cite{kansas}(or see \cite{fenni}) as

\begin{equation} \label{}
f(z) = 1 + \frac{2}{\pi^2} \Big(\ln(1+\sqrt{z})/(1-\sqrt{z})\Big)^2 .
\end{equation}

Note that we may simplify by substituting $z$ for $\sqrt{z}$, since $||f(z)||_{H^2} = ||f(z^2)||_{H^2}$ (or see Proposition \ref{littleguy} below). The coefficients of this series can be found by squaring the series for

\begin{equation} \label{}
\ln(1+z)/(1-z) = 2(z+\frac{z^3}{3} + \frac{z^5}{5} + \ldots ) = \sum_{n=0}^{\ff} \frac{2z^{2n+1}}{2n+1} .
\end{equation}

We obtain

\begin{equation} \label{HHpar}
f(z^2) = 1 + \frac{2}{\pi^2} \Big(\ln(1+z)/(1-z)\Big)^2 = 1 + \frac{8}{\pi^2} \sum_{n=0}^\ff \Big(\sum_{j=0}^{n} \frac{1}{(2j+1)(2(n-j)+1)}\Big) z^{2n+2} .
\end{equation}

We apply Proposition \ref{bigguy} with $a=1$, using \rrr{pareq} and \rrr{HHpar}, to obtain

\bea \label{}
&& \nn 1^2 \! + \! \Big(1 \! \times \! \frac{1}{3} + \frac{1}{3} \! \times \! 1 \Big)^2 \! + \! \Big(1 \! \times \! \frac{1}{5} + \frac{1}{3} \! \times \! \frac{1}{3} + \frac{1}{5} \! \times \! 1 \Big)^2 \! + \! \Big(1 \! \times \! \frac{1}{7} + \frac{1}{3} \! \times \! \frac{1}{5} + \frac{1}{5} \! \times \! \frac{1}{3} + \frac{1}{7} \! \times \! 1 \Big)^2 + \dots
\\ \nn && \hspace{2.5cm} = \sum_{n=0}^\ff \Big(\sum_{j=0}^{n} \frac{1}{(2j+1)(2(n-j)+1)}\Big)^2 = \frac{\pi^4}{32} .
\eea

\vski

\ccases{hyper} Let $V_p$ be the hyperbolic region defined by

\begin{equation} \label{}
\Big\{ \Big(\frac{\sin^2(\frac{\pi}{2}p)x + \cos^2(\frac{\pi}{2}p)}{\cos(\frac{\pi}{2}p)}\Big)^2 - \sin^2(\frac{\pi}{2}p)y^2 > 1, \quad x>0\Big\},
\end{equation}

where $p$ is a parameter with $0<p<1/2$. The reader may check that

\begin{equation} \label{}
\GG_V(x+iy) = \frac{2\cos^2(\frac{\pi}{2}p)}{\sin^2(\frac{\pi}{2}p)(\cos^2(\frac{\pi}{2}p)-\sin^2(\frac{\pi}{2}p))}\Big(\Big(\frac{\sin^2(\frac{\pi}{2}p)x + \cos^2(\frac{\pi}{2}p)}{\cos(\frac{\pi}{2}p)}\Big)^2 - \sin^2(\frac{\pi}{2}p)y^2 - 1\Big).
\end{equation}

The conformal map for $\DD$ to $V_k$ sending 0 to 1 (which is a focus of the hyperbola) is given in \cite{kansas} as

\begin{equation} \label{}
f(z) = \frac{1}{\sin^2(\frac{\pi}{2}p)} \cosh\Big(p \log \frac{1+\sqrt{z}}{1-\sqrt{z}}\Big) - \cot^2(\frac{\pi}{2}p).
\end{equation}

As in Example \ref{para}, we may simplify by substituting $z$ for $\sqrt{z}$. We calculate

\begin{equation} \label{}
\begin{split}
f(z^2) & = \frac{1}{2\sin^2(\frac{\pi}{2}p)}\Big(\Big(\frac{1+z}{1-z}\Big)^{p} + \Big(\frac{1-z}{1+z}\Big)^{p}\Big)- \cot^2(\frac{\pi}{2}p) \\
& = \frac{1}{2\sin^2(\frac{\pi}{2}p)}\Big(\frac{(1+z)^{2p} + (1-z)^{2p}}{(1-z^2)^{p}} \Big)- \cot^2(\frac{\pi}{2}p) \\
& = \frac{1}{2\sin^2(\frac{\pi}{2}p)}\Big(\sum_{j=0}^\ff \binom{2p}{j} z^j + \sum_{j=0}^\ff \binom{2p}{j} (-z)^j\Big) \Big(\sum_{j=0}^{\ff}\binom{p+j-1}{j}z^{2j}\Big) - \cot^2(\frac{\pi}{2}p) \\
& = \frac{1}{2\sin^2(\frac{\pi}{2}p)}\Big(2 \sum_{j=0}^\ff \binom{2p}{2j} z^{2j}) \Big(\sum_{j=0}^{\ff}\binom{p+j-1}{j}z^{2j}\Big) - \cot^2(\frac{\pi}{2}p) \\
& = \frac{1}{\sin^2(\frac{\pi}{2}p)}\sum_{n=0}^\ff \Big(\sum_{j=0}^{n}\binom{2p}{2j} \binom{p+n-j-1}{n-j}\Big)z^{2n} - \cot^2(\frac{\pi}{2}p),
\end{split}
\end{equation}

where we have applied the generalized Binomial Theorem in the following forms, where $\binom{a}{b} = \frac{a(a-1)\ldots (a-b+1)}{b!}$ is the binomial coefficient:

\begin{equation} \label{bino}
(1+z)^{2p} = \sum_{j=0}^\ff \binom{2p}{j} z^j, \quad (1-z^2)^{-p} = \sum_{j=0}^{\ff}\binom{p+j-1}{j}z^{2j}.
\end{equation}

Applying Proposition \ref{bigguy} with $a=1$ and simplifying gives

\begin{equation} \label{free}
\sum_{n=1}^\ff \Big(\sum_{j=0}^{n}\binom{2p}{2j} \binom{p+n-j-1}{n-j}\Big)^2 = \frac{2\sin^4(\frac{\pi}{2}p)}{\cos^2(\frac{\pi}{2}p)-\sin^2(\frac{\pi}{2}p)} = \frac{2\sin^4(\frac{\pi}{2}p)}{\cos(\pi p)}
\end{equation}

It is clear that this approaches $\ff$ as $p \nearrow \frac{1}{2}$, which means that the sum on the left side of \rrr{free} diverges for $1/2 \leq p \leq 1$. This follows also from the results in \cite{burk}, which show that no function satisfying (i)-(iii) can exist on $V_p$ if $p \in [1/2,1)$ (see Theorem 3.1 of \cite{burk} and the ensuing application). We note that two values of $p$ which give simple expressions on the right side of \rrr{free} are $p=1/3$, which yields

\begin{equation} \label{}
\sum_{n=1}^\ff \Big(\sum_{j=0}^{n}\binom{2/3}{2j} \binom{1/3+n-j-1}{n-j}\Big)^2 = \frac{1}{4},
\end{equation}
and $p=1/4$, which yields

\begin{equation} \label{}
\sum_{n=1}^\ff \Big(\sum_{j=0}^{n}\binom{1/2}{2j} \binom{1/4+n-j-1}{n-j}\Big)^2 = \frac{1}{4}(3\sqrt{2}-4).
\end{equation}

\vski

\ccases{hyper2} Let $V_\th$ be the region between the two branches of the hyperbola given by $(\frac{x}{\sin \th})^2 - (\frac{y}{\cos \th})^2 = 1$, where $\th \in (0,\pi/4)$. That is, $V=\{(\frac{x}{\sin \th})^2 < (\frac{y}{\cos \th})^2 + 1\}$. It is not hard to see that

\begin{equation} \label{}
\GG_V(x+iy) = \frac{2\sin^2 \th \cos^2 \th}{\cos^2 \th - \sin^2 \th}\Big(1+(\frac{y}{\cos \th})^2 - (\frac{x}{\sin \th})^2\Big) .
\end{equation}

It is shown in \cite{fenni} that, for fixed $\th$, $f(\th,z)=\sin(\frac{4\th}{\pi}\tan^{-1} z)$ is a conformal map from $\DD$ to $V_\th$ mapping $0$ to $0$. If we expand $f(\th,z) = \sum_{n=1}^{\ff} a_n(\th)z^n$, then Proposition \ref{bigguy} applied with $a=0$ shows that

\begin{equation} \label{truth}
\sum_{n=1}^{\ff} a_n(\th)^2 = \frac{2\sin^2 \th \cos^2 \th}{\cos^2 \th - \sin^2 \th} = \frac{\sin^2 2\th}{2 \cos 2\th}.
\end{equation}

A value for $\th$ where we may find an explicit expression for $a_n(\th)$ is $\th= \frac{\pi}{8}$. The half-angle formula for sine leads to

\begin{equation} \label{}
f(\pi/8,z) = \sin(\frac{1}{2}\tan^{-1} z) = \sqrt{\frac{1-\cos (\tan^{-1} z)}{2}} =\sqrt{\frac{1-\frac{1}{\sqrt{1+z^2}}}{2}} = \frac{\sqrt{\sqrt{1+z^2}-1}}{\sqrt{2}(1+z^2)^{1/4}}.
\end{equation}

Applying the Binomial Theorem in the form of \rrr{bino} gives

\begin{equation} \label{}
(1+z^2)^{-1/4} = \sum_{j=0}^{\ff}\binom{j-3/4}{j}(-1)^jz^{2j}.
\end{equation}

It is also true that

\begin{equation} \label{}
\sqrt{\sqrt{1+z^2}-1} = \frac{z}{\sqrt{2}} \; \; {}_2F_1(1/4,3/4;3/2;-z^2) = \frac{1}{\sqrt{2}}\sum_{j=0}^{\ff} \frac{(-1)^jC(2j)z^{2j+1}}{16^j},
\end{equation}

where $C(k) = \frac{1}{k+1} \binom{2k}{k}$ is the $k$-th Catalan number (see \cite{oeis} for details). We obtain

\begin{equation} \label{}
\sin(\frac{1}{2}\tan^{-1} z) = \frac{1}{2}\sum_{n=0}^{\ff} (-1)^n \Big(\sum_{j=0}^{n} \binom{j-3/4}{j}\frac{C(2(n-j))}{16^{n-j}}\Big) z^{2n+1}.
\end{equation}

Thus, \rrr{truth} gives the identity

\begin{equation} \label{}
\sum_{n=0}^{\ff} \Big(\sum_{j=0}^{n} \binom{j-3/4}{j}\frac{C(2(n-j))}{16^{n-j}}\Big)^2 = \sqrt{2}.
\end{equation}

Similarly to Case \ref{hyper}, it can be shown that no function satisfying (i)-(iii) can exist on $V_\th$ for $\th \in [\frac{\pi}{4},\frac{\pi}{2})$.

\vski

\ccases{ell} Let $V_\xi$ be the elliptical region given by $\{\frac{x^2}{\cosh^2 \xi} + \frac{y^2}{\sinh^2 \xi} < 1\}$, where $\xi$ is a parameter greater than 0. The ellipse determining $V_\xi$ has foci at $-1$ and $1$. We may take

\begin{equation} \label{}
\GG_{V_\xi}(x+iy) = \frac{-2\sinh^2 \xi \cosh^2 \xi}{\sinh^2 \xi + \cosh^2 \xi}\Big(\frac{x^2}{\cosh^2 \xi} + \frac{y^2}{\sinh^2 \xi} - 1\Big).
\end{equation}

Let $K(z,t)$ be the normal elliptic integral of the first kind (following the normalization employed in \cite{fenni}), that is

\begin{equation} \label{}
K(z,t) = \int_{0}^{z} \frac{dx}{\sqrt{(1-x^2)(1-t^2x^2)}}
\end{equation}

for $0<t<1$. Set

\begin{equation} \label{}
\mu(t) = \frac{\pi K(1,\sqrt{1-t^2})}{2K(1,t)}.
\end{equation}

Suppose $t$ is chosen such that $\mu(t)=2\xi$, and let $f_{a,t}(z)$ be the conformal map from $\DD$ onto $V_\xi$ mapping 0 to $a$ with $f_{a,t}'(0)>0$. It is known that

\begin{equation} \label{}
f_{0,t}(z) = \sin\Big(\frac{\pi}{2K(1,t)}K(z/\sqrt{t},t)\Big);
\end{equation}

see \cite{fenni} for details. Let

\begin{equation} \label{}
f_{0,t}(z) = \sum_{n=1}^{\ff}a_n(t) z^n.
\end{equation}

It is shown in \cite{fenni} that the coefficients $a_n(\xi)$ are positive for odd $n$ and 0 for even $n$. Proposition \ref{bigguy} then shows that upon setting $A_n(t) = a_{2n+1}(t)$ we have

\begin{equation} \label{}
\sum_{n=0}^{\ff}A_n(t)^2 = \frac{2\sinh^2 \xi \cosh^2 \xi}{\sinh^2 \xi + \cosh^2 \xi} = \frac{\sinh^2 2\xi }{2\cosh 2\xi} = \frac{\sinh^2 \mu(t) }{2\cosh \mu(t)}.
\end{equation}

Furthermore, it is shown in \cite{fenni} that

\begin{equation} \label{}
f_{-1,2\xi}(z^2) = 2 (f_{0,\xi}(z))^2 - 1.
\end{equation}

Thus,

\begin{equation} \label{}
f_{-1,2\xi}(z^2) = -1 + 2 \sum_{n=0}^{\ff} \Big(\sum_{j=0}^{n}A_j A_{n-j}\Big)z^{2n+2}.
\end{equation}

Applying Theorem \ref{bigguy} with $a = -1$ yields

\begin{equation} \label{}
\sum_{n=0}^{\ff} \Big(\sum_{j=0}^{n}A_j A_{n-j}\Big)^2 = \frac{\sinh^2 2\xi \cosh^2 2\xi}{2(\sinh^2 2\xi+\cosh^2 2\xi)}\Big(\frac{\cosh^2 2\xi-1}{\cosh^2 2\xi}\Big) = \frac{\sinh^4 \mu(t)}{2\cosh 2\mu(t)}.
\end{equation}

Further identities involving the values of $A_n(t)$ are possible using higher order Chebyshev polynomials in order to construct maps from $\DD$ to $V_{m\xi}$; see \cite[p. 333]{fenni}. It is also shown in \cite{fenni} that the values $A_n(t)$ can be characterized as the unique solution to the recurrence relation

\begin{gather}
(2n+2)(2n+3)A_{n+1}(t) - \Big((t+\frac{1}{t})(2n+1)^2-\frac{\pi^2}{4tK(1,t)^2}\Big)A_n(t) + 2n(2n-1)A_{n-1}(t) = 0. \\
\nn A_{-1}=0, A_0 = \frac{\pi}{2 \sqrt{t} K(1,t)}.
\end{gather}

These values can take a relatively simple form when $t$ is chosen so that $K(1,t)$ is a rational multiple of $\pi$; for instance, if $\om$ is such that $K(1,\om) = \pi$, then we have

\begin{gather} \label{}
A_0(\om) = \frac{1}{2\sqrt{\om}}, \quad A_1(\om) = \frac{3+4\om^2}{48\om^{3/2}}, \quad A_2(\om) = \frac{105 + 56 \om^2 + 144 \om^4}{3840\om^{5/2}}, \\ \nn A_3(\om) = \frac{10395 + 4524 \om^2 + 4496 \om^4 + 14400 \om^6}{645120 \om^{7/2}}, \quad \ldots
\end{gather}

Unfortunately, however, it seems difficult to find an explicit formula for the coefficients $A_n(t)$.

\section{Concluding remarks}

\begin{remark}
Cases \ref{para} through \ref{ell} contain the formulas for the expectation of the first time that planar Brownian motion hits each of the conics.
\end{remark}

\begin{remark}
An examination of the proof of Proposition \ref{bigguy} will show that the injectivity of $f_a$ is not crucial to the argument. Instead, the important features possessed by conformal maps are analyticity and the mapping of the boundary of $\DD$ to the boundary of $V$. With this in mind, let us define a function $f$ from $\DD$ into $V$ to be {\it B-proper} if $f$ is analytic and $\lim_{r \nearrow 1} f(re^{i\th})$ exists and lies in $\dd V$ for almost every $\th$ in $[0,2\pi)$. A similar definition was made in \cite{meecp} (the definitions agree if $||f||_{H^2}<\ff$). The proof of Proposition \ref{bigguy} therefore applies to give us the following proposition.

\begin{prop} \label{littleguy}
Let $f_a(z) = a + \sum_{n=1}^{\ff}a_n z^n$ be a B-proper map from $\DD$ into $V$. Suppose that $u_V$ is a continuous function on $V$ satisfying (i$'$)-(iii$'$) from Proposition \ref{bigguy}. Then $u_V(a) = \sum_{n=1}^{\ff} |a_n|^2$.
\end{prop}

It is difficult to find a really good example to which we can apply this result, but we may give a simple one as follows. Let $f(z) = e^{tan^{-1}z}$. $f$ is not conformal, but does map $\DD$ analytically onto the annulus $\{e^{\frac{-\pi}{4}} < |z| < e^{\frac{\pi}{4}}\}$, taking 0 to 1. $f$ is B-proper, since $f$ extends continuously to map $\{z=e^{i\th}; -\frac{\pi}{2}<\th <\frac{\pi}{2}\}$ to $\{ |z| = e^{\frac{\pi}{4}}\}$ and $\{z=e^{i\th}; \frac{\pi}{2}<\th <\frac{3\pi}{2}\}$ to $\{ |z| = e^{\frac{-\pi}{4}}\}$. We may therefore apply Proposition \ref{littleguy}. The function $\ln |z|$ is harmonic on $V$, so it may be directly verified that we may take

\begin{equation} \label{}
\GG(a) = \Big(\frac{e^{\pi/2} - e^{-\pi/2}}{\pi/2}\Big)\ln |a| + \frac{e^{\pi/2} + e^{-\pi/2}}{2} - |a|^2 .
\end{equation}

The power series for $f(z)$ begins $e^{\tan^{-1} z} = 1 + z + \frac{z^2}{2} - \frac{z^3}{6} - \frac{7z^4}{24} + \frac{z^5}{24} + \ldots$. We conclude from Proposition \ref{littleguy} that

\begin{equation} \label{rant}
1 + \frac{1}{2^2} + \frac{1}{6^2} + \frac{7^2}{24^2} + \frac{1}{24^2} + \ldots = \frac{e^{\pi/2} + e^{-\pi/2}}{2} - 1.
\end{equation}

It is unfortunately not easy to find the general expression for the terms in this series. However, we may characterize the terms to be those determined uniquely by the following recurrence relation:

\begin{equation} \label{}
a_n - (n+1)a_{n+1} - (n-1)a_{n-1} = 0, \qquad a_1 = 1.
\end{equation}

This relation comes from equating the coefficients in the power series obtained in the identity

\begin{equation} \label{}
(1+z^2) \frac{d}{dz} e^{tan^{-1}z} = e^{tan^{-1}z}.
\end{equation}
\end{remark}

\begin{remark}
Suppose we have a domain $V$, a formula for $\GG_V(a)$, and knowledge of a conformal map $f_{a_o}$ from $\DD$ to $V$ sending 0 to $a_o$. Then we obtain not only an identity concerning the power series of $f_{a_o}$ but in fact an entire family of identities, one for each point in $V$, as setting $w=f_{a_o}^{-1}(a)$ the map $f(\frac{z+w}{1+\bar w z})$ is a conformal map sending $\DD$ to $V$ and $0$ to $a$. However, in practice the power series of the functions $f(\frac{z+w}{1+\bar w z})$ are difficult to explicitly compute, limiting the utility of this observation. We do, however, obtain a bit of intuition on how the method works. We are realizing a specific sum as one in a family of sums, indexed by a complex parameter $a$. In certain cases it is easier to evaluate the entire family of sums in one operation, using Propositions \ref{bigguy1} or \ref{bigguy}, than it is to evaluate the specific case we are interested in. We then recover the special case by choosing the correct value for $a$.
\end{remark}

\begin{remark}
The recent paper \cite{coffey} contains several more examples to which this method applies.
\end{remark}

\section{Acknowledgements}

I would like to thank George Markowsky and Mark Coffey for useful conversations, as well as an anonymous referee for helpful comments. I am also grateful for support from Australian Research Council Grant DP0988483.


\def\noopsort#1{} \def\printfirst#1#2{#1} \def\singleletter#1{#1}
   \def\switchargs#1#2{#2#1} \def\bibsameauth{\leavevmode\vrule height .1ex
   depth 0pt width 2.3em\relax\,}
\makeatletter \renewcommand{\@biblabel}[1]{\hfill#1.}\makeatother

\bibliographystyle{alpha}
\bibliography{CABMbib}

\providecommand{\bysame}{\leavevmode\hbox to3em{\hrulefill}\thinspace}
\providecommand{\MR}{\relax\ifhmode\unskip\space\fi MR }
\providecommand{\MRhref}[2]{%
  \href{http://www.ams.org/mathscinet-getitem?mr=#1}{#2}
}
\providecommand{\href}[2]{#2}
\begin{thebibliography}{10}

\bibitem{oeis}
\emph{The {O}n-{L}ine {E}ncyclopedia of {I}nteger {S}equences},
  http://oeis.org/A048990.

\bibitem{ahl}
L.V. Ahlfors, \emph{Complex analysis: {A}n introduction to the theory of
  analytic functions of one complex variable}, McGraw-Hill, 1966.

\bibitem{ala}
A.~Alabert, M.~Farr\'{e}, and R.~Roy, \emph{Exit times from equilateral
  triangles}, Applied Mathematics and Optimization \textbf{49} (2004), no.~1,
  43--53.

\bibitem{drum}
R.~Ba{\~n}uelos and T.~Carroll, \emph{Brownian motion and the fundamental
  frequency of a drum}, Duke Mathematical Journal \textbf{75} (1994), no.~3,
  575--602.

\bibitem{burk}
D.L. Burkholder, \emph{Exit times of {B}rownian motion, harmonic majorization,
  and {H}ardy spaces}, Advances in Mathematics \textbf{26} (1977), no.~2,
  182--205.

\bibitem{coffey}
M.~Coffey, \emph{Expected exit times of {B}rownian motion from planar domains:
  {C}omplements to a paper of {M}arkowsky}, arXiv:1203.5142v1.

\bibitem{dritref}
T.A. Driscoll and L.N. Trefethen, \emph{Schwarz-{C}hristoffel mapping},
  Cambridge University Press, Cambridge, UK, 2002.

\bibitem{schaumpde}
P.~DuChateau and D.W. Zachmann, \emph{Schaum's outline of theory and problems
  of partial differential equations}, Schaum's Outline Series, 1986.

\bibitem{hardwood}
G.H. Hardy and JE~Littlewood, \emph{A maximal theorem with function-theoretic
  applications}, Acta Mathematica \textbf{54} (1930), no.~1, 81--116.

\bibitem{helm}
K.~Helmes, S.~R\"ohl, and R.H. Stockbridge, \emph{Computing moments of the exit
  time distribution for {M}arkov processes by linear programming}, Operations
  Research (2001), 516--530.

\bibitem{fenni}
S.~Kanas and T.~Sugawa, \emph{On conformal representations of the interior of
  an ellipse}, Annales Academi\ae { }Scientiarum Fennic\ae { } \textbf{31}
  (2006), 329--348.

\bibitem{kansas}
S.~Kanas and A.~Wisniowska, \emph{Conic regions and k-uniform convexity},
  Journal of Computational and Applied Mathematics \textbf{105} (1999),
  no.~1-2, 327--336.

\bibitem{knight}
F.B. Knight, \emph{Essentials of {B}rownian motion and diffusion}, American
  Mathematical Society, 1981.

\bibitem{woepf}
W.~Koepf, \emph{Schwarz-{C}hristoffel mappings: Symbolic computation of mapping
  functions for symmetric polygonal domains}, Proceedings of the Workshop on
  Functional-Analytic Methods in Complex Analysis and Applications to Partial
  Differential Equations (1995), 293--305.

\bibitem{meecp}
G.~Markowsky, \emph{On the expected exit time of planar {B}rownian motion from
  simply connected domains}, Electronic Communications in Probability
  \textbf{16} (2011), 652--663.

\bibitem{rud}
W.~Rudin, \emph{Real and complex analysis}, Tata McGraw-Hill, 2006.

\end{thebibliography}

\end{document}